# Further Remarks on the Sampled-Data Feedback Stabilization Problem

John Tsinias and Dionysios Theodosis

*Abstract*— The paper deals with the problem of the sampled data feedback stabilization for autonomous nonlinear systems. The corresponding results extend those obtained in [3] by the same authors. The sufficient conditions we establish are based on the existence of discontinuous control Lyapunov functions and the corresponding results are applicable to a class of nonlinear affine in the control systems.

## I. INTRODUCTION

We consider autonomous systems of the general form:

$$\dot{x} = f(x,u), (x,u) \in \mathbb{R}^n \times \mathbb{R}^m, \quad f(0,0) = 0 \quad (1)$$

where $f: \mathbb{R}^n \times \mathbb{R}^m \to \mathbb{R}^n$ is Lipschitz continuous. For every initial $x_0 \in \mathbb{R}^n$, measurable and locally essentially bounded control $u:[s,t_{\max}) \to \mathbb{R}^m$, we denote by $\pi(\cdot) = \pi(\cdot,s,x_0,u)$ the corresponding trajectory of (1), that satisfies $\pi(s,s,x_0,u) = x_0$, where $t_{\max} = t_{\max}(s,x_0,u)$ is the corresponding maximal existance time of the trajectory. We say that system (1) is Semi-Globally Asymptotically Stabilizable by Sampled-Data Feedback (SDF-SGAS), if for every constant $R > 0$ and for any given partition of times $T_1 := 0 < T_2 < T_3 < \ldots < T_\nu < \ldots$, with $T_\nu \to \infty$, there exist a neighborhood $\Omega_R$ of zero with

$$B[0,R;\mathbb{R}^n] := \{x \in \mathbb{R}^n : |x| \le R\} \subset \Omega_R,$$

(where $|x|$ denotes the Euclidean norm of the vector $x \in \mathbb{R}^n$) and a map $k: \mathbb{R}^+ \times \Omega_R \to \mathbb{R}^m$ such that for any $x \in \Omega_R$ the map $k(\cdot,x): \mathbb{R}^+ \to \mathbb{R}^m$ is measurable and locally essentially bounded, zero is stable with respect to the sampled-data closed-loop system $\dot{\pi} = f(\pi,k(t,\pi(T_i)))$, $t \in [T_i, T_{i+1})$, $i = 1,2,\ldots$ with $\lim_{t \to \infty} \pi(t) = 0$, $\forall \pi(0) \in \Omega_R$. The definition above is adopted in [3,5] and constitutes a time-varying version of both concepts of asymptotic controllability and sampled-data stabilization (see [1,2]).

J. Tsinias, is with the Department of Mathematics, National Technical University of Athens, Athens 15780, Greece (e-mail:, jtsin@central.ntua.gr). D. Theodosis is with the Dynamic Systems and Simulation Laboratory, Technical University of Crete, Chania, Greece (e-mail: dtheodosis@dssl.tuc.gr).

The following property modifies the usual concept of "control Lyapunov function" adopted in [1,2,3] and other related works.

**Property 1:** For system (1) assume that there exist $a_1, a_2, a \in K$ ($a_1, a_2, a$ are continuous, strictly increasing with $a_1(0) = a_2(0) = a(0) = 0$) and $V: \mathbb{R}^n \to \mathbb{R}^+$ (being in general discontinuous) such that

$$a_1(|x|) \le V(x) \le a_2(|x|), \quad x \in \mathbb{R}^n \quad (2)$$

in such a way that for every $R > 0$, $x \in B[0,R;\mathbb{R}^n], x \ne 0$, there exist constants $\sigma_x, L_x, M_x > 0$ such that for every $\varepsilon \in (0,\sigma_x]$, there exists an input

$$u_{\varepsilon,x} := u_{\varepsilon,x}(\cdot):[0,\varepsilon] \to B[0,M_x;\mathbb{R}^m] \quad (3)$$

with

$$V(\pi(\varepsilon,0,\bar{x},u_{\varepsilon,x})) < V(\bar{x}) - L_x; \quad (4a)$$

$$V(\pi(t,0,\bar{x},u_{\varepsilon,x})) \le a(V(\bar{x})), \quad (4b)$$

$$\forall t \in [0,\varepsilon], \bar{x} \text{ near } x$$

The following lemma generalizes the result in [3, Proposition 2] establishing SDF-SGAS for (1) under the stronger hypothesis that $V$ is continuous.

**Lemma 1:** Property 1 implies SDF-SGAS for system (1).

**Remark 1:** The establishment of the result above is based on the same procedure applied in the proof of [3, Proposition 2] pluss certain modifications; we note that the following facts play a central role for the proof of Lemma 1.

**Fact 1:** Due to our assumption that, in addition to (4), condition (3) is fulfilled for certain $M_x \in (0,+\infty)$, the corresponding trajectory involved in (4) also satisfies:

$$|\pi(t,0,\bar{x},u_{\varepsilon,x}) - \bar{x}| \le \varepsilon C_x \quad (5)$$

$$\forall t \in [0,\varepsilon], \bar{x} \text{ near } x, \varepsilon > 0 \text{ near zero}$$

for certain $C_x > 0$.

**Fact 2:** The pair of conditions (3) and (4) are equivalent to the following property: For every nonempty compact set $S \subset \mathbb{R}^n \setminus \{0\}$ there exist constants $\sigma, L > 0$, such that for every $R > 0$, $x \in S \cap B[0,R;\mathbb{R}^n]$, $x \ne 0$ and $\varepsilon \in (0,\sigma]$ there exists a constant $M_x > 0$ and input $u = u_{\varepsilon,x}$ satisfying (3) and such that

$$V(\pi(\varepsilon,0,x,u_{\varepsilon,x})) < V(x) - L; \quad (6a)$$

$$V(\pi(t,0,x,u_{\varepsilon,x})) \le 2a(V(x)); \quad \forall t \in [0,\varepsilon] \quad (6b)$$

**Remark 2:** If $V$ is lower semicontinuous, then (6) is equivalent to the weaker assumption that for the specific $x$ both (4a,b) hold with $L_x = 0$. In particular, instead of (4a,b), we may assume that for any sufficiently small $\varepsilon \in (0, \sigma_x]$ there exists a control $u_{\varepsilon,x} : [0, \varepsilon] \to B[0, M_x; \mathbb{R}^m]$ with

$$V(\pi(\varepsilon, 0, x, u_{\varepsilon,x})) < V(x); \qquad (7a)$$

$$V(\pi(t, 0, x, u_{\varepsilon,x})) \leq 2a(V(x)), \quad \forall t \in [0, \varepsilon]. \qquad (7b)$$

It turns out that *continuity* of $V$ implies equivalence between (4a,b), (6a,b) and (7a,b).

Lemma 2 of Section II is the main result of this article and establishes that, under existence of an appropriate family of *continuous* Lyapunov-like functions $V_i$, with $\cup \text{dom} V_i = \mathbb{R}^n$, a control Lyapunov function W can be found, being *discontinuous* and satisfying all conditions of Property 1, including (4) with $V$ instead of $W$. The latter according to Lemma 1 guarantees solvability of SDF-SGAS for (1). We use the result of Lemma 2 to derive sufficient conditions for the solvability of SDF-SGAS for a class of affine in the control systems. (Section III: Proposition 1, Corollary 1 and Proposition 2).

## II. THE MAIN RESULT

We present the main result of this article based on the existence of a family of continuous control Lyapunov function establishing sufficient conditions for SDF-SGAS systems systems (1). Its proof is based on the result of Lemma 1.

**Lemma 2 (Main Result):** For system (1) assume that there exist $\omega_1, \omega_2 \in K$ and $a \in K$, such that for every constant $R \in \mathbb{R}^+$ there exist open regions

$$A_i \subset \mathbb{R}^n \setminus \{0\}, i = 1, 2, \ldots. \qquad (8)$$

with

$$A_i \cap A_j = \emptyset, \forall i \neq j \qquad (9a)$$

$$B[0, R; \mathbb{R}^n] \subset cl\left(\cup A_i\right) \qquad (9b)$$

and continuous mappings

$$V_i : \mathbb{R}^n \to \mathbb{R}^+, V_i(0) = 0, i = 1, 2, \ldots \qquad (10)$$

with

$$\omega_1(|x|) \leq V_i(x) \leq \omega_2(|x|), i = 1, 2, \ldots, x \in A_i \qquad (11)$$

and in such a way that for every $i = 1, 2, \ldots$, nonzero $x$ belonging to $cl(A_i)$, and sufficiently small $\varepsilon > 0$, there exist a constant $M_i > 0$ and a control

$$u^i_{\varepsilon,x} : [0, \varepsilon] \to B[0, M_i; \mathbb{R}^m] \qquad (12)$$

with

$$V_i(\pi(\varepsilon, 0, x, u^i_{\varepsilon,x})) < V_i(x); \qquad (13a)$$

$$V_i(\pi(t, 0, x, u^i_{\varepsilon,x})) \leq a(V_i(x)), \forall t \in [0, \varepsilon] \qquad (13b)$$

Then, system (1) satisfies all conditions of Lemma 1 for some appropriate Upper Semicontinuous (USC) Lyapunov function $W : \mathbb{R}^n \to \mathbb{R}^+$ satisfying (2)-(4) with $V := W$, therefore, according to Lemma 1, the system (1) is SDF-SGAS.

**Remark 3:** For the case where there exist a finite number of open sets $A_i \subset \mathbb{R}^n \setminus \{0\}$, $i = 1, 2, \ldots, N$, $N \in \mathbb{N}$ such that (9.a,b) hold and associated with continuous functions $V_i : \mathbb{R}^n \to \mathbb{R}^+$ as in statement of Proposition 2, then *positive definiteness* of each $V_i : A_i \to \mathbb{R}^+, i, \ldots, N$ implies (11).

**Proof of Lemma 2:** Without any loss of generality, we may assume that all regions $A_i$ are bounded and, due to (9) and (10), for every $x \neq 0$ there exists an integer $k$ and finite number of indices $i_1, i_2, i_3, \ldots, i_k$ for which $x \in cl(A_{i_1} \cup A_{i_2} \ldots \cup A_{i_k})$ and $x \notin cl(A_j)$, $\forall j \neq i_1, i_2, \ldots, i_k$. Consider for each $i = 1, 2, \ldots$, a constant

$$c_i = c_i(A_i) > 0 \qquad (14)$$

and functions $a_1, a_2 \in K$ such that

$$a_1(|x|) \leq \omega_1(|x|) + \min\{c_i \text{ for those } i \text{ for which } x \in cl(A_i)\} \qquad (15a)$$

$$\omega_2(|x|) + \max\{c_i \text{ for those } i \text{ for which } x \in cl(A_i)\} \leq a_2(|x|) \qquad (15b)$$

$$a(V_i(x)) + c_i < 2a(V_i(x) + c_i); \quad \forall x \in cl(A_i) \qquad (15c)$$

where $a$ is defined in (13b), and in such a way that, if we define

$$W_i := V_i + c_i, x \in A_i \qquad (16)$$

the following holds

$$W_i(x) \neq W_j(x), \forall x \in cl(A_i) \cap cl(A_j), x \neq 0 \qquad (17)$$

Define

$$W(x) \begin{cases} := W_i(x), & x \in A_i \\ := \max\{W_j(x); j = i_1, i_2, \ldots, i_k \in \mathbb{N}\} \text{ for} \\ (0 \neq) x \in \cup \partial A_j \partial A_{i_1} \cup \partial A_{i_2} \ldots \cup \partial A_j : \\ x \notin \cap \partial A_j; j \neq i_1, i_2, \ldots, i_k \text{ for certain } k = k(x) \in \mathbb{N} \\ := 0 & x = 0 \end{cases} \qquad (18)$$

and note that (16) and (17) guarantee existence of such a $k$. The function $W$ is USC and, due to (8)- (18), conditions (2)-(4) hold with $V := W$, $a := 2a$ and $a_1$, $a_2$ instead of $\omega_1$, $\omega_2$, respectively. More specifically, (2) is an immediate consequence of (8), (11), (14)-(18). In order to show (4), we consider two cases:

**Case 1:** $x \in A_i$, $x \neq 0$, for some integer $i = i(x)$.

It follows by invoking our hypotheses that for every sufficiently small $\varepsilon > 0$, there exist a constant $M_i > 0$ and an input $u^i_{\varepsilon,x}:[0,\varepsilon] \to B[0,M_i;\mathbb{R}^m]$ satisfying both (13a,b) and (12). The latter implies:

$$|\pi(t,0,x,u_{\varepsilon,x}) - x| \leq \varepsilon C_i, \forall t \in [0,\varepsilon], \varepsilon > 0 \quad (19)$$

for certain $C_i > 0$. It follows that for sufficiently small $\varepsilon > 0$ and $x \in A_i$ the trajectory $\pi(t,0,x,u^i_{\varepsilon,x})$ remains inside $A_i$ for $t$ near zero. By recalling (12)-(18) and selecting sufficiently small $\varepsilon > 0$, we have:

$$x \in A_i \Rightarrow$$
$$W(\pi(\varepsilon,0,x,u^i_{\varepsilon,x})) = W_i(\pi(\varepsilon,0,x,u^i_{\varepsilon,x})) \quad (20)$$
$$= V_i(\pi(\varepsilon,0,x,u^i_{\varepsilon,x})) + c_i < V_i(x) + c_i = W(x)$$

and simultaneously:

$$x \in A_i \Rightarrow$$
$$W(\pi(t,0,x,u^i_{\varepsilon,x})) = W_i(\pi(t,0,x,u^i_{\varepsilon,x}))$$
$$= V_i(\pi(t,0,x,u^i_{\varepsilon,x})) + c_i \leq a(V_i(x)) + c_i \quad (21)$$
$$\leq 2a(V_i(x) + c_i) = 2a(W(x)), t \in [0,\varepsilon]$$

for appropriate $u^i_{\varepsilon,x}$ satisfying (12). Then, for the specific $x$ above, the desired (4a,b) are consequence of (19)-(21) and continuity of $V_i$. Particularly, (4) is valid with $V := W, a := 3a$, sufficiently small $L_x > 0$ and $a_i, i = 1,2$ as above.

**Case 2:**

$$(x \neq 0), x \in \partial A_{i_1} \cap \partial A_{i_2} ... \cap \partial A_{i_k}, j = i_1, i_2, ..., i_k \in \mathbb{N}:$$
$$x \notin \bigcap \partial A_j; j \neq i_1, i_2, ..., i_k \text{ for certain } k = k(x) \in \mathbb{N}$$

For the specific $x$ and $i$ as above and by taking into account (15)-(18), and (22) we may define the integer

$$I^i_x := \max \left\{ p : W_p(x) = \max \left\{ \begin{array}{l} W_j(x); j = i_1, i_2, ..., i_k \in \mathbb{N}: \\ x \in \partial A_{i_1} \cap \partial A_{i_2} ... \cap \partial A_{i_k} : x \notin \partial A_j, \\ j \neq i_1, i_2, ..., i_k \end{array} \right\} \right\} \quad (22)$$

Obviously by (17), (18), (19), and (22), we have

$$W(x) = W_{I^i_x}(x) \quad (23a)$$

$$W_{I^i_x}(x) \geq W_j(x), \forall j \neq i_1, i_2, ..., i_k \quad (23b)$$

and, due to continuity of $V_i$,

$$W_{I^i_x}(y) > W_j(x), \forall j \neq I^i_x, \quad (24a)$$

$$W_{I^i_x}(y) \geq W_j(x), \forall j \neq i_1, i_2, ..., i_k,$$
$$I^i_x = I^i_y, \quad (24b)$$
$$\forall y : |y - x| \leq \xi$$

for appropriately small $\xi := \xi(x) > 0$ The latter, in conjunction with (18) implies

$$W(y) = W_{I^i_x}(y), \forall |y - x| \leq \xi \quad (25)$$

We next recall our assumptions (12) and (13) in order, for the specific $i$ and $x$ above, to find a constant $M_i > 0$ and a control satisfying (3) and (4). Indeed, by (12) and (13) we may pick sufficiently small $\varepsilon > 0$ and determine a control

$$u^i_{\varepsilon,x} := u^{I^i_x}_{\varepsilon,x}[0,\varepsilon] \to B[0,M_i;\mathbb{R}^m] \quad (26)$$

in such a way that the trajectory $x(t,0,x,u^i_{\varepsilon,x})$ satisfies (12) and (13). Notice that, due to (24b), $u^{I^i_x}_{\varepsilon,x} = u^{I^i_y}_{\varepsilon,y}$ and $W(\pi(\varepsilon,0,y,u^i_{\varepsilon,x})) = W(y)$ for y near x. It follows from (12), (13), (15c), (18), (23)-(26) that for every $i_1, i_2, ..., i_k$ as above we have:

$$0 \neq x \in \partial A_i \Rightarrow$$
$$W(\pi(\varepsilon,0,y,u^i_{\varepsilon,y})) = W_{I^i_y}(\pi(\varepsilon,0,y,u^i_{\varepsilon,x}))$$
$$\leq V_{I^i_x}(\pi(\varepsilon,0,y,u^i_{\varepsilon,x})) + c_{I^i_x}$$
$$< V_{I^i_x}(y) + c_{I^i_x}$$
$$= \max \left\{ \begin{array}{l} W_j(y), j = i_1, i_2, ..., i_k \in \mathbb{N}; y \in \partial A_{i_1} \cap \partial A_{i_2} ... \cap \partial A_{i_k}, \\ y \notin \bigcap \partial A_j; j \neq i_1, i_2, ..., i_k \in \partial A_k \end{array} \right\}$$
$$= W(y), y \text{ near } x$$
$$(27)$$

and simultaneously

$$W(\pi(t,0,y,u^i_{\varepsilon,y})) = W_{I^i_y}(\pi(\varepsilon,0,y,u^i_{\varepsilon,y}))$$
$$\leq V_{I^i_x}(\pi(t,0,y,u^i_{\varepsilon,x})) + c_{I^i_x}$$
$$\leq a\left(V_{I^i_x}(y)\right) + c_{I^i_x}$$
$$\leq 2a\left( \max \left\{ \begin{array}{l} W_j(y), j = i_1, i_2, ..., i_k \in \mathbb{N}; y \in \partial A_{i_1} \cap \partial A_{i_2} ... \cap \partial A_{i_k}, \\ y \notin \bigcap_k \partial A_j; j \neq i_1, i_2, ..., i_k \in \partial A_k \end{array} \right\} \right)$$
$$= 2a(W(x)), y \text{ near } x$$
$$(28)$$

For the specific $x$ above and appropriate small $L_x > 0$ the desired (4a,b) are consequences of (26)-(28), continuity of $V_i$ on the region $A_i$. This completes the proof.

## III. APPLICATIONS

We apply the result of Lemma 2 to some interesting cases of nonlinear systems. Consider the affine in the control single input system:

$$\dot{x} = f(x) + ug(x), x \in \mathbb{R}^n$$
$$f(0) = 0 \quad (29)$$

where $f, g$ are smooth $(C^\infty)$. The main result below combines the results of Lemma 2 of present work and the algebraic result established in [3]. We first provide some standard appropriate notations: By $Lie\{f,g\}$ we denote the Lie algebra generated by $\{f,g\}$. Let $L_1 = span\{f,g\}$ and $L_{i+1} = span\{[X,Y], X \in L_i, Y \in L_1\}$, $i=1,2,...$ and for every nonzero $\Delta \in Lie\{f,g\}$ define:

$$order_{\{f,g\}}\Delta \begin{cases} :=1, \text{ if } \Delta \in L_1 \setminus \{0\} \\ :=k>1, \text{ if } \Delta = \Delta_1 + \Delta_2, \text{with } \Delta_1 \in L_k \setminus \{0\} \\ \qquad \text{and } \Delta_2 \in span\{\cup_{i=1}^{i=k-1} L_i\} \end{cases}$$

**Proposition 1:** For (29) we assume that there exist $\alpha_1, \alpha_2, \alpha \in K$, open nonempty regions $A_i \subset \mathbb{R}^n \setminus \{0\}$ with $cl(\cup A_i) = \mathbb{R}^n$, and smooth mappings $V_i : \mathbb{R}^n \to \mathbb{R}^+$ such that $a_1(|x|) \le V_i(x) \le a_2(|x|)$, $i=1,2..$, and for each $i=1,2,...$ the following hold: For every $x \in A_i$ either

$$(gV_i)(x) \ne 0, \qquad (30a)$$

or one of the following holds: Either

$$(gV_i)(x) = 0 \Rightarrow (fV_i)(x) < 0 \qquad (30b)$$

or there exists an integer $N = N(i,x) \ge 1$ such that

$$(gV_i)(x) = 0, (f^j V_i)(x) = 0, \; j = 1,2,...,N;$$

$$(\Delta_{i_1} \Delta_{i_2} \ldots \Delta_{i_k} V_i)(x) = 0$$
$$\forall \Delta_{i_1}, \Delta_{i_2}, \ldots, \Delta_{i_k} \in Lie\{f,g\} \setminus \{g\} \qquad (30c)$$
$$\text{with } \sum_{p=1}^{k} order_{\{f,g\}} \Delta_{i_p} \le N$$

where $(f^j V_i)(x) := f(f^{j-1}V_i)(x)$, $i = 2,3,...,$ $(f^1 V_i)(x) := (fV_i)(x)$ and in such a way that one of the following properties hold:

P1 $\qquad (f^{N+1}V_i)(x) < 0$

P2 $N$ is odd and

$$([[...[[f,g],g],...,g],g]V_i)(x) \ne 0$$
$$\underbrace{\qquad}_{N \text{ times}}$$

P3 $N$ is even and

$$([[...[[f,g],g],...,g],g]V_i)(x) < 0$$
$$\underbrace{\qquad}_{N \text{ times}}$$

P4 $N$ is an arbitrary positive integer with

$$(f^{N+1}V_i)(x) = 0,$$

$$([[...[[g,f],f],...,f],f]V_i)(x) \ne 0$$
$$\underbrace{\qquad}_{N \text{ times}}$$

Then, conditions of Lemma 2 are satisfied,-with same $a(s) := 2s$ in all cases- hence, (1) is SDF-SGAS.

**Proof(Outline):** The proof of Proposition 1 is based on the fact that, if the algebraic conditions imposed above hold, then, according to [3, Proposition 2], for any $R>0$ all conditions (9)-(13) of Lemma 2 are satisfied with $V, A_i$ as above with same $a(s) := 2s$ in all cases, hence, system (29) is SDF-SGAS. •

The following result constitutes a special case of Proposition 1 for systems:

$$\begin{pmatrix} \dot{x} \\ \dot{y} \end{pmatrix} = f + ug := \begin{pmatrix} F(x,y) \\ 0 \end{pmatrix} + u \begin{pmatrix} 0 \\ 1 \end{pmatrix} \qquad (31)$$

where $F : \mathbb{R}^n \times \mathbb{R} \to \mathbb{R}^n$ is $C^1$ with $F(0,0) = 0$ and generalize [4, Theorem 1] (see also [5]) concerning the stabilization problem by means of a "feedback integrator".

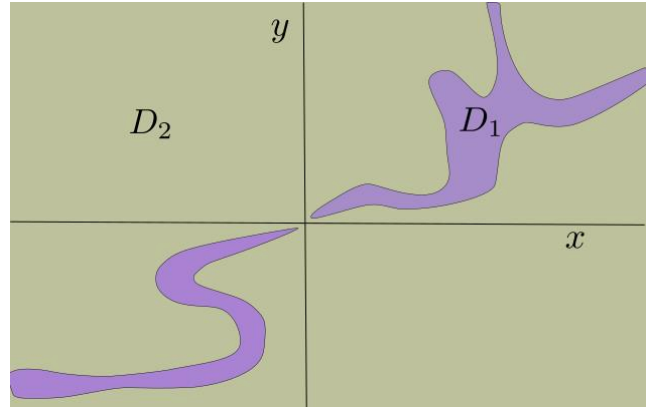

**Figure**: Illustration of the sets $D_1$ and $D_2$ in Corollary 1.

We make the following assumption

**Corollary 1:** For system (31), assume that there exists a $C^1$ positive definite function $V : \mathbb{R}^n \to \mathbb{R}^+$ with:

$$V(0) = 0, \; V(x) > 0, \; \forall x \ne 0 \qquad (32)$$

a $C^1$ function $W : \mathbb{R}^n \times \mathbb{R} \to \mathbb{R}^+$ with

$$W(0,0) = 0 \qquad (33)$$

and open sets $D_1, D_2 \subset \mathbb{R}^{n+1} \setminus \{0\}$ such that

$$D_1 \cap D_2 = \varnothing \qquad (34a)$$

$$D_1 \cup D_2 = \mathbb{R}^{n+1} / \{0\} \qquad (34b)$$

satisfying the following properties:

$$(x,y) \ne 0, \; (x,y) \in cl(D_1) \Rightarrow x \ne 0 \qquad (35)$$

$$(x,y) \ne 0, \; (x,y) \in cl(D_1) \Rightarrow$$

$$\text{either} \quad DV(x)F(x,y) < 0 \qquad (36a)$$

$$\text{or} \quad \begin{cases} DV(x)F(x,y) = 0; \\ DV(x)\dfrac{\partial F}{\partial y}(x,y) \ne 0 \end{cases} \qquad (36b)$$

$$(0, y) \in cl(D_2); W(0, y) = 0 \Rightarrow y = 0 \quad (37)$$

$$0 \neq (x, y) \in D_2 \Rightarrow \frac{\partial W}{\partial y}(x, y) \neq 0 \quad (38)$$

Then, under the previous properties, system (31) satisfies the assumptions of Proposition 1, and therefore is SDF-SGAS.

**Proof**: Define:

$$V_1 := V \quad (39)$$

We have

$$DV_1(x)(f + ug)\big|_{(x,y) \in cl(D_1) \setminus \{0\}} = DV_1(x)F(x, y)\big|_{(x,y) \in cl(D_1) \setminus \{0\}}$$

which in conjunction with (36) and (38) implies that, either

$$DV_1(x)(f + ug) < 0, \text{ for all } u \quad (40a)$$

or

$$\begin{cases} DV_1(x)(f + ug) = 0 \\ DV_1(x)[g, f] = DV(x)\frac{\partial}{\partial y}F(x, y) < 0 \end{cases}, \text{ for all } u \quad (40b)$$

From (32), (35) and (38), it follows that $V_1$ satisfies (11) on the region $D_1 \setminus \{0\}$ for certain $\omega_1$ and $\omega_2$. From (40) it also follows that the pair $(A_1, V_1)$, $A_1 := D_1$ satisfies (30) of Proposition 1 for the system (31). Define next:

$$V_2(x, y) := V_1(x) + W(x, y) \quad (41)$$

From (41) we have:

$$DV_2(x, y)(f + ug)\big|_{(x,y) \in cl(D_2) \setminus \{0\}}$$
$$= \left(\frac{\partial V_1}{\partial x} + \frac{\partial W}{\partial x}, \frac{\partial W}{\partial y}\right)\begin{pmatrix} F(x, y) \\ u \end{pmatrix}\bigg|_{(x,y) \in cl(D_2) \setminus \{0\}} \quad (42)$$
$$= \left(\frac{\partial V_1}{\partial x}F(x, y) + \frac{\partial W}{\partial x}F(x, y) + u\frac{\partial W}{\partial y}\right)\bigg|_{(x,y) \in cl(D_2) \setminus \{0\}}$$

From (37),(38) and (42) it also follows that the pair $(A_2, V_2)$, $A_2 := D_2$ satisfies (30a,b) plus P2 of Proposition 1 - with N=1 - for the system (31) and further (11) holds for certain $\omega_1$ and $\omega_2$. The latter, in conjunction with (34) and Remark 3, implies that, all desired properties of Proposition 1 are satisfied with N=1, therefore system (31) is SDF-SGAS. •

**Remark 4:** The sufficient condition of Corollary 1 can be weakened by replacing $DV(x)\frac{\partial F}{\partial y}(x, y) = DV(x)[f, g](x, y) \neq 0$ in (36a) by the general assumption that for the function $V$ above there exists an integer $N \geq 1$ such that (30c) plus one the conditions P1, P2, P3, P4 are fulfiled, where $f$ and $g$ are defined in (31).

Finally, we deal with the SDF-SGAS problem for systems (1) with dynamics

$$f(x, u) = A(x)x + B(x)u \\ x \in \mathbb{R}^n, u \in \mathbb{R}^m \quad (43)$$

where $A : \mathbb{R}^n \to \mathbb{R}^n$, $B : \mathbb{R}^n \to \mathbb{R}^{n \times m}$ are Lipschitz and assume that for every nonzero $x \in \mathbb{R}^n$ there exists a map $F = F(x) \in \mathbb{R}^{m \times n}$ such that

$$(A + BF)(x) \text{ is Hurwitz} \quad (44)$$

equivalently, we assume that for every nonzero $x \in \mathbb{R}^n$ the pair $(A(x), B(x))$ is stabilizable. As a consequence of Lemma 2 we get:

**Proposition 2:** Under (44), system (43) is SDF-SGAS.

**Proof:** According to our assumption (44), for every $\xi \in \mathbb{R}^n$ there exist a symmetric $P = P(\xi) \in \mathbb{R}^{n \times n}$ which is positive definite, namely,

$$x'P(\xi)x > 0, \ \forall x \neq 0 \quad (45a)$$

and further

$$x'P(\xi)(A(\xi) + B(\xi)F(\xi))x \leq -k(\xi)|x|^2, \ \forall x \in \mathbb{R}^n \quad (45b)$$

for certain $k = k(\xi) > 0$. In addition, for any $R > 0$ a pair of constants $c, C > 0$ can be found such that

$$c \leq |P(\xi)| \leq C, \ \forall \xi \in B[0, R; \mathbb{R}^n] \quad (45c)$$

Define

$$V_\xi(x) := \frac{1}{2}x'P(\xi)x \quad (46a)$$

$$u_\xi(x) := F(\xi)x \quad (46b)$$

and consider the trajectory $x(\cdot) = \pi(\cdot, 0, \xi, u_\xi)$, $x(0) = \xi$, of the system

$$\dot{x} = A(x)x + B(x)u_\xi \\ = (A(x) + B(x)F(\xi))x \quad (47)$$

Then, by evaluating the time-derivative $\dot{V}_\xi$ of $V_\xi = \frac{1}{2}x'(t)P(\xi)x(t)$ along the trajectory of (47). We find from (45b) and (46):

$$\dot{V}_\xi = x'(t)P(\xi)\big((A(x(t)) - A(\xi)) + A(\xi) \\ + (B(x(t)) - B(\xi))F(\xi) + B(\xi)F(\xi)\big)x(t) \\ = x'(t)P(\xi)(A(\xi) + B(\xi)F(\xi))x(t) \\ + x'(t)P(\xi)(A(x(t)) - A(\xi))x(t) \quad (48) \\ + x'(t)P(\xi)(B(x(t)) - B(\xi))F(\xi)x(t) \\ \leq |x(t)|^2 \big(-k(\xi) + |P(\xi)||A(x(t)) - A(\xi)| \\ + |P(\xi)||B(x(t)) - B(\xi)||F(\xi)|\big)$$

which implies:

$$V_\xi(x(t)) \leq V_\xi(x(0)) - \frac{1}{2}k(\xi)|x(t)|^2,$$

$\forall t \geq 0$ near zero in such a way that $x(t)$ lies in a neighborhood $A_\xi$ of $\xi$ (49)

By exploiting (45c), (48) and (49) for every constant $R \in \mathbb{R}^+$ we can determine $\xi_1, \xi_2, \ldots \in \mathbb{R}^n$ and open neighborhoods $A_i = A_{\xi_i}$, $i = 1, 2 \ldots$ of $\xi_i$ and appropriate positive definite $C^1$ functions

$$V_i(x) = \frac{1}{2}x'P(\xi_i)x, \quad x \in A_i$$

such that all conditions (8)-(13) are satisfied with $a = 1$ and for certain $\omega_1, \omega_2 \in K$. We conclude that (43) satisfies all conditions of Lemma 2, therefore is SDF-SGAS. •

### REFERENCES

[1] F.H. Clarke, Y.S. Ledyaev, E.D. Sontag and A.I. Subbotin, ``Asymptotic controllability implies feedback stabilization'', *IEEE Transactions on Automatic Control*, 42, 10, 1394-1407, 1997.

[2] E.D. Sontag, *Mathematical control theory*, $2^{nd}$ edn., Springer, Berlin, Heidelberg, New York, 1998.

[3] J.Tsinias and D. Theodosis, "Sufficient Lie Algebraic Conditions for Sampled-Data Feedback Stabilizability of Affine in the Control Nonlinear Systems, *IEEE Transactions on Automatic Control*, 61, 5, 1334-1339, 2016.

[4] J. Tsinias, "An extension of Artstein's, Theorem on Stabilization by Ordinary Feedback Integrators", *Systems & Control Letters*, 20, 141-148, 1993.

[5] J. Tsinias, "Remarks on Asymptotic Controllability and Sampled-Data Stabilization for Autonomous Systems*", IEEE Transactions on Automatic Control, 55, 3, 721-726, 2010.*